\documentclass[12pt]{article}
\def\Ima{\mathrm{Im}\, }
\def\R{\mathbb{R}}
\def\C{\mathbb{C}}
\usepackage{amssymb,amsfonts,amsthm,amsmath}
\usepackage{graphicx, verbatim}
\usepackage{xcolor}
\title{Iossif Ostrovskii's work on entire functions}
\author{Alexandre Eremenko and Mikhail Sodin}
\begin{document}
\maketitle
\begin{abstract}
The theory of entire functions and its applications were at the center of Ostrovskii's research interests throughout his entire career. He made lasting contributions to several aspects of this theory, and many of his works had a significant influence on subsequent research. In this note, we address some of these works.

\end{abstract}

\noindent
{\bf 1. Growth of entire functions and their value distribution by arguments}
\vspace{.1in}

We start with Ostrovskii's earliest results, included in his PhD thesis, in which
he extended and strengthened seminal results
of Krein and Edrei. He published announcements of these results along with sketches of the proofs in a series of Doklady notes (1957-1960). A detailed exposition appeared in the Izvestiya paper~\cite{Ostr-Izvestiya}. Most of them, in somewhat stronger versions, were included in the Goldberg--Ostrovskii treatise~\cite[Sections~VI.2, VI.3]{GO}.

\vspace{.1in}
In 1947, Krein introduced the class of entire functions
$f$ such that $1/f$ is represented by an absolutely convergent series of simple fractions
\[
\frac{1}{f(z)} = a + \frac{b}{z} + \sum_{\lambda\in\Lambda}  c_\lambda 
\Bigl( \frac1{z-\lambda} + \frac1\lambda \Bigr),
\]
where $a$, $b$, and $c_\lambda$ are real, $\Lambda\subset\R$, and
$ \sum_{\Lambda} |c_\lambda|/\lambda^2 <\infty $. We will denote this class by $\mathsf K$.
Krein proved that 

\vspace{.1in}
{\em Functions of the class
$\sf K$ have finite exponential type}.

\vspace{.1in}
\noindent 
Moreover,  Krein showed that {\em functions from $\sf K$
have a bounded type in the  upper and lower half-planes},
that is, in each of the half-planes they are represented as a
quotient of two bounded analytic functions.
The latter  property yields that functions from $K$ belong to the Cartwright
class, which consists of entire functions $f$ of exponential type with a growth bound on the real axis:
\[
\int_\R \frac{\log^+ |f(t)|}{1+t^2}\, dt < \infty.
\]
Functions of the Cartwright class possess a regularity of growth and of zero distribution described in books by
Levin~\cite[Lectures~14--17]{Levin} and Koosis~\cite{Koosis-Cart}. 

It is also important to mention that Krein proved his  theorem  
under an assumption weaker than $\Lambda\subset\R$, namely that the zeros of $f$ satisfy the Blaschke condition in the
upper and lower half-planes
\begin{equation}\label{B}
\sum_{\lambda\in\Lambda} \bigl| \Ima (1/\lambda)\bigr| < \infty.
\end{equation}

Edrei's theorem~\cite[Theorem~1]{Ed} deals with meromorphic functions with three 
radially distributed values, which is a seemingly different class of functions. We will quote here only 
its special case  pertaining to entire functions in which case only two valued are needed (the third one is $\infty$):

\vspace{.1in}
{\em If $f$ is an entire function such that for some distinct  $a, b\in\C$ all but finitely many 
solutions to the equations 
\[
f(z)=a, \qquad f(z)=b
\]
lie in a finite union of rays 
\[
D(\alpha_1\, \ldots \, \alpha_N) = \bigcup_{1\le j \le N } \{z=re^{ i\alpha_j}\colon 0\le r<\infty\}, \quad 
0\le \alpha_1 <\,  \ldots  \, < \alpha_N  < 2\pi,
\] 
then $f$ has a finite order
of growth
\begin{equation}\label{eq:rho}
\rho_f \le \max_{1\le j \le N}\, \frac{\pi}{\alpha_{j+1}-\alpha_j},
\end{equation}
where $\alpha_{N+1}=2\pi + \alpha_1$}.

\vspace{.1in}
\noindent 
In particular, if an entire function $f$ has real $\pm 1$-points, then its order of growth does not
exceed one. We will return  to this class of entire functions in Section~3. 

It should be mentioned that the bound~\eqref{eq:rho} under the
a priori assumption that $f$ has
a finite order of growth is implicitly contained in the work of Bieberbach~\cite{Bieber}.

\vspace{.1in}
A remarkable joint feature of these theorems of Krein and Edrei is that both of them do not have 
a priori assumptions on the growth of the entire function. Their similarity, communicated to
Ostrovskii by Goldberg, allowed Ostrovskii to significantly extend both of them.
Here, we will bring only special cases of his results. The interested reader can find more general
versions in~\cite[Section~VI.2]{GO}.

Following Ostrovskii, we say that $a$-points $z_k = r_k e^{ i\theta_k}$ 
of an entire function $f$ are {\em very close}
to the system of rays $D(\alpha_1\, \ldots \, \alpha_N)$, if they satisfy the Blaschke condition
in each of the complementary angles to these rays. That is, the $N$  series
\[
\sum_{ \substack {r_k \ge 1, \\ 
\theta_k \in (\alpha_j, \alpha_{j+1}) } } 
r_k^{-\frac{\pi}{\alpha_{j+1}-\alpha_j}}\, 
\sin \left[\frac{\pi}{\alpha_{j+1}-\alpha_j}\, (\theta_k-\alpha_j) \right],
\qquad j=1, 2, \ldots , N,
\]
converge (in the case, when $N=2$ and $\alpha_1=0$, $\alpha_2=\pi$, this coincides with
condition~\eqref{B}). Then, his extension of Krein's theorem states the following:

\vspace{.1in}
{\em 
Suppose that $f$ is an entire function such that $1/f$ is 
represented by an absolutely convergent series of simple fractions with
poles very close to a finite union of rays $D(\alpha_1\, \ldots \, \alpha_N)$.

Then $f$ has a finite order of growth satisfying~\eqref{eq:rho} 
and has a bounded type in each of the angles complementary to $D(\alpha_1\, \ldots \, \alpha_N)$.
}

\vspace{.1in}
\noindent
We note that the fact that $f$ has a bounded type in each of the complementary angles yields
regularity of growth and zero distribution similar to the ones which are possessed by entire functions of the 
Cartwright class, see~\cite[Sections~VI.2, VI.3]{GO}. 

\vspace{.1in}
\noindent
Ostrovskii's strengthening of Edrei's theorem is also quite natural:

\vspace{.1in}
{\em 
Suppose that $a$, $b$ are distinct values in $\C$, and that $D(\alpha_1\, \ldots \, \alpha_N)$ 
is a finite union of rays. Let $f$ be an entire function with $a$- and $b$-points lying very close to  
 $D(\alpha_1\, \ldots \, \alpha_N)$.

Then $f$  has a finite order of growth satisfying~\eqref{eq:rho} and
has a bounded type in each of the complementary angles.
}

\vspace{.1in}
Both results have meromorphic counterparts also proven by Ostrovskii, see~\cite[Section~VI.2]{GO}
and~\cite{Ostr-Izvestiya} for the case of meromorphic functions in the unit disk.
For the reader acquainted with basic notions of the Nevanlinna theory, we mention that 
if the function $f$ is meromorphic, then the same conclusion holds if $a$- and $b$- points  
lie very close to  $D(\alpha_1\, \ldots \, \alpha_N)$, $a, b \in\bar \C$, $a\ne b$, and a value $c\in\bar\C$, $c\ne a, b$,
is Nevanlinna deficient. The proofs are based on a version of Nevanlinna theory for meromorphic functions in
an angle. 

Edrei and Fuchs~\cite{EF-1962} considered meromorphic functions with all but finitely many 
$a$- and $b$-points ($a \ne b\in\C$) lying on $N$ disjoint rectifiable curves 
$C_j =\{z=re^{ i\alpha_j (r)}, r_0 \le r <\infty  \}$, 
where $\alpha_1(r)<\, \ldots \,  < \alpha_N(r) < \alpha_1(r)+2\pi = \alpha_{N+1}(r)$.
They assumed that, for some $B>0$,  the length 
of the portion of each curve lying in any annulus $\{r_1\le |z|\le r_2\}$ is bounded by $B (r_2-r_1)$. They showed that
if $f$ has relatively few poles (more precisely they assumed that $\infty$ is a Nevanlinna deficient value of $f$),
then $f$ has a finite order of growth, which does not exceed
\[
\max_{1\le j \le N}\, \lim_{r\to\infty}\, \frac{9\pi B^2}{\alpha_{j+1}(r)-\alpha_j(r)}\,. 
\]
If the curves are radial lines emanating from the origin, then $B=1$, and their bound is $9$ times worse than 
the optimal one~\eqref{eq:rho}.

\vspace{.1in}
The interested reader will find other results (including Ostrovskii's ones) on the value-distribution of meromorphic
functions in~\cite[Chapter~VI]{GO} as well as in the survey~\cite[\S~10]{GLO}.

\vspace{.1in}
Concluding this section, it is worth noting that the Goldberg--Ostrovskii treatise~\cite{GO}, written in the late 1960s, 
remains one of the primary sources in one-dimensional Nevanlinna theory. Its English translation includes a brief 
survey of results obtained after 1970. A less direct but still important impact of Ostrovskii on value-distribution theory 
may be found in works of Petrenko, who completed his Ph.D. thesis under Ostrovskii's supervision. 
Many of Petrenko's findings were summarized in his books~\cite{Petrenko1, Petrenko2}.

\vspace{.2in}
\noindent
{\bf 2. Conjectures of P\'olya and Wiman}
\vspace{.1in}

The class of all real polynomials with all zeros real is
closed under differentiation. A simple proof of this statement uses Rolle's
theorem and degree counting. The second part of the argument, degree counting,
breaks down when one tries to extend this result to entire functions. The entire
function $ze^{z^2}$ has only real zeros while the derivative $(1+2z^2)e^{z^2}$
has two non-real zeros. This justifies the definition of the
Laguerre--P\'olya class of entire functions: $f\in \sf{LP}$ if $f$ is a limit
(uniform on compact subsets of the plane) of real polynomials with all zeros
real. It immediately follows that this class is closed under differentiation. Laguerre and P\'olya obtained 
a remarkable parametric description of this class: $f\in \sf{LP}$ if and only if  $f(z)=e^{-az^2}  g(z)$,
where $a\ge 0$ and $g$ is a real entire function (i.e. having real Taylor coefficients)
of genus one with real zeros. That is,
$$f(z)=z^me^{-az^2+bz+c}\prod_k \Bigl( 1-\frac{z}{z_k} \Bigr) e^{z/z_k},$$
where $b,c$ and $z_k$ are real,  $\sum_k |z_k|^{-2} <\infty$, 
and $a\geq 0$. 
This class plays a central role in the ``algebraic theory of entire functions'' launched 
by Laguerre, P\'olya, and Schur, and has a variety of applications.

In 1914, P\'olya~\cite{P2} stated the following conjecture:

\vspace{.1in}
\noindent 
P1. {\em Suppose $f$ is a real entire function with all derivatives having only real zeros. Then
$f\in\sf{LP}$}.

\vspace{.1in}\noindent
In~\cite{P2} (a more detailed exposition appeared in~\cite{P2a}) P\'olya proved this 
for real entire functions of the form $f=Pe^Q$, where $P$ and $Q$ are polynomials. 
In~\cite{P2a}  P\'olya gave a stronger version of this conjecture which includes non-real functions:

\vspace{.1in}
\noindent 
P2. {\em Suppose $f$ is an entire function with all derivatives having only real zeros. Then
either $f$ has the form
\[
f(z) = ce^{az} \quad {\rm or} \quad f(z) = c(e^{ ibz} - e^{ id}),
\]
with complex $c$ and $a$, and real $b$ and $d$, or $f\in {\sf LP}$}.

\vspace{.1in}
In~\cite{A, A-C.R.} {\AA}lander claimed the proofs of Conjectures~P1 and~P2 for entire functions 
of finite order. Subsequent researchers expressed doubts about the justification of certain steps in his proofs, 
see~\cite[p. 228]{HW-1}. 

According to {\AA}lander~\cite{Ad},  Wiman conjectured a stronger version of Conjecture~P1
in his Uppsala lectures of 1911:

\vspace{.1in}
\noindent 
W1. {\em Suppose $f$ is a real entire function with $f$ and $f''$ having only real zeros. Then
$f\in\sf{LP}$}.

\vspace{.1in}
Furthermore, according to {\AA}lander, Wiman also gave a quantitative version of  
this conjecture for entire functions of finite order of growth. To state it we define the class ${\sf V}_{2p}$ 
consisting of all entire functions of the form
$f(z)=e^{-az^{2p+2}}g(z)$,
where $a\ge 0$, and  $g$ is a real entire function  of genus not exceeding  $2p+1$ having only real zeros.
Then, we set ${\sf W}_0 = {\sf V}_0$, and ${\sf W}_{2p}={\sf V}_{2p}\backslash {\sf V}_{2p-2}$, so that  
${\sf LP}={\sf W}_0$, and the set of all entire functions of finite order with only real zeros is a disjoint union
of classes ${\sf W}_{2p}$ for all $p\geq 0$.

\vspace{.1in}
\noindent 
W2. 
{\em For $f\in {\sf W}_{2p}$, $f''$ has at least $2p$ non-real zeros. }

\vspace{.1in} In other words, if  $f$ is a real entire functions with only real zeros and $f''$ has 
at most $2p$ non-real zeros, then $f\in {\sf V}_{2p}$. In particular, for $p=0$, this coincides with 
Conjecture~W1 for a special case of entire functions of finite order. 
In~\cite{Ad},  {\AA}lander \cite{Ad} proved Conjecture W2 for functions of genera 2 and 3, and in his 
later paper \cite{A0} for functions of genera not exceeding 5.

To see the difficulty of these conjectures, consider the case $f=e^Q$ with a real polynomial $Q$
and let $P=Q'$.
Conjecture~W1 yields that {\em for every real polynomial $P$ of degree $n>1$,
$P^2+P'$ has non-real
zeros}. Moreover, if ${\rm deg} P = n$, then the genus of $f$ equals $n+1$ and
$f\in {\sf W}_{2p}$, where $p$ is the least integer greater or equal to $n/2$.
Then, by Conjecture~W2, 
the polynomial $P^2+P'$ of degree $2n$ should have at least $2p$ non-real zeros. This is
known for polynomials $P$ having only real zeros, see P\'olya--Szeg\H{o}~\cite[Chapter~5, 182]{PS},
but a general, even a non-quantitative, question remained open for about 60 years. But more on that later.

\vspace{.1in}
In the paper~\cite{PFourier} P\'olya was concerned with ``a very bold argument'' by which Fourier tried to prove
that all zeros of the Bessel function $J_0$ are real and positive.  Fourier applied a rule proven only for polynomials,
to a transcendental entire function, while Cauchy and Poisson expressed doubts about Fourier's reasoning.
In~\cite{PFourier} P\'olya proved a theorem which justifies Fourier's argument and stated hypothetical theorems
which would allow one to broadly apply Fourier's algebraic argument. 
We will refer to these two hypothetical theorems, which we will state shortly,  as the Fourier--P\'olya conjectures.

Let $f$ be a real entire function. We assume for simplicity that neither $f$ nor any of its derivatives $f'$, $f''$,  \ldots \,
has multiple zeros. We say that a point $\xi\in\R$ is {\em de Gua-critical} if, for some $n\in\mathbb N$, 
\[
f^{(n)}(\xi) = 0, \quad f^{(n-1)}(\xi)f^{(n+1)}(\xi) >0\,.
\]
Counting de Gua-critical points, we take them into account as many times as this occurs for different values of $n$.
That is, we count all positive minima and negative maxima of $f$ and all its derivatives. 
The de Gua rule (1741) states that the number of complex conjugate
pairs of zeros of a real polynomial equals the number
of its de~Gua-critical points. This rule is not difficult to justify using, for instance, the idea of the Budan-Fourier rule, while Fourier applied it to a real transcendental function, justifying this by
the fact that this function is a limit of real polynomials.

\vspace{.1in}
\noindent
FP1. {\em A real entire function of genus $0$ has just as many de Gua critical points as pairs of non-real zeros. 
}

\vspace{.1in}
\noindent To state the second conjecture, we extend the Laguerre-P\'olya class allowing the
functions to have finitely many non-real zeros, and denote by  ${\sf LP}^*$ the class of entire functions of the form $f=pg$, where $p$ is a real polynomial, and $g\in{\sf LP}$.

\vspace{.1in}
\noindent
FP2. {\em An entire function of class ${\sf LP}^*$
has just as many de Gua critical points as pairs of non-real zeros. 
}

\vspace{.1in}
\noindent
In the same paper P\'olya shows that the second conjecture is equivalent to another one, to which we refer 
as to a PW conjecture, since, according to {\AA}lander~\cite{A-Wiman}, Wiman stated it to him in 1916.

\vspace{.1in}
\noindent
PW. {\em  If $f\in{\sf LP^*}$, then its derivatives from a certain one onward, have no non-real zeros at all
(that is, belong to ${\sf LP}$).
 }

\vspace{.1in} 
Special cases of this conjecture were proven by {\AA}lander~\cite{A-Wiman}, P\'olya~\cite{PFourier, P-1937}, 
and Wiman~\cite{W1, W2}, so it was established for entire functions of order at most $4/3$, and for entire functions of
the form $f(z)=e^{-\alpha z^2} g(z)$, where $\alpha\ge 0$ and $g$ has genus $0$.  

In an enthralling survey~\cite{P4}, P\'olya briefly returned to this circle of problems and
stated yet another conjecture, which concludes this list:

\vspace{.1in}
\noindent
P3. {\em If the order of a real entire function $f(z)$ is greater than 2, and $f(z)$
has only a finite number of non-real zeros, then the number of non-real
zeros of $f^{(n)}$ tends to infinity as $n\to\infty$.}

\vspace{.1in}
All these conjectures were completely settled much later, the most difficult ones only in the 21st Century, 
through the combined efforts of many mathematicians. 
The first substantial progress was made in the paper \cite{LO} by Levin and Ostrovskii. Most of the 
subsequent work on the conjectures built upon this paper, where Levin and Ostrovskii introduced several crucial 
ideas to the subject. The first one was the connection between the logarithmic derivatives
of entire functions with all zeros real and the class of meromorphic functions
with positive real part in the upper half-plane. 

\vspace{.1in}
\noindent
{\bf Levin-Ostrovskii representation:} {\em Let $f$ be a real entire function
with all zeros real. Then
$$\frac{f'}{f}=Pg,$$
where $P$ is real entire and $g$ is a real meromorphic
function with non-negative imaginary part in the upper half-plane.
If $f$ is of finite order, then $P$ is a polynomial.}
\vspace{.1in}

This lemma gives a control of the behavior of $f'/f$ since $g$ has
nice properties. One obtains especially good control for functions of finite
order, in which case the representation can be made more precise
\cite{HW-2, BE}.

Then, following Edrei~\cite{Ed}, Levin and Ostrovski consider the function
$F=f/f'$ and argued as follows. Suppose that $ff''$ has no zeros in the upper half-plane. Then
$F$ omits $0$ and $F'$ omits $1$ in the upper
half-plane. Indeed,
$$F'-1=\frac{f'^2-ff''}{(f')^2}-1=-\frac{ff''}{(f')^2}.$$
This observation immediately points at  Hayman's generalization of Picard's theorem~\cite[\S~3.3]{H}:

\vspace{.1in}
{\em If a meromorphic function in the plane omits one finite value
and its derivative omits a non-zero value, then this function is constant}.

\vspace{.1in}
\noindent Applied to the function $F=f/f'$, Hayman's theorem yields that if $ff''$ has no zeros   
in $\C$, then $f$ is the exponential function $e^{Az+B}$.
There is a general philosophy, usually called Bloch's Principle, that if a 
condition imposed on a function meromorphic in the plane implies that this
function is constant, then the same condition imposed on a function in an
arbitrary domain must imply some universal estimate of this function.

Hayman's proof relies on Nevanlinna theory, and to apply it to a function
in a half-plane required an appropriate generalization of Nevanlinna theory.
Versions of Nevanlinna characteristics adapted for
functions in a half-plane were introduced by Nevanlinna, and by Levin and Tsuji.
The formulation with Levin--Tsuji characteristics is especially suitable, since for them a full analogue
of the main technical tool of Nevanlinna theory (the Lemma on the logarithmic
derivative) holds. Using both Nevanlinna characteristics for a half-plane
and Levin--Tsuji characteristics, Levin and Ostrovskii proved the following:

\vspace{.1in}
{\em If $f$ is a real entire function such that $ff''$ has only real zeros,
then}
\begin{equation}
\label{1}
\log\log M(r,f)=O(r\log r).
\end{equation}
Here $M(r,f)=\max\{|f(z)|:|z|\leq r\}$.
A weaker estimate of the same type under the stronger condition
that $ff'f''$ has only real zeros was earlier obtained by Edrei~\cite[Theorem~3]{Ed}.
It is worth mentioning that Levin and Ostrovskii proved \eqref{1} requiring that zeros
of $f$ are real while zeros of $f''$ satisfy Blaschke's condition~\eqref{B}
in the upper and lower half-planes. Shen~\cite{Shen} proved that~\eqref{1}
holds provided that $f$ is real, and zeros of $f$ and $f''$ satisfy~\eqref{B}.

\vspace{.1in}
Below we list the main milestones on the subsequent way to the complete proof of
the conjectures of P\'olya and Wiman.

\vspace{.1in}\noindent
In 1971 Hellerstein and Yang \cite{HY} extended the Levin--Ostrovskii
theorem to higher derivatives, that is, if $ff^{(k)}$ has only real zeros for some $k\geq 2$, then (\ref{1}) holds.

\vspace{.1in}\noindent
In 1977, Hellerstein and Williamson~\cite{HW-1, HW-2} proved that if $ff'f''$ has only real zeros
then $f\in LP$. This proves P\'olya's conjecture P1, even in a stronger version, and yields that
the Laguerre-P\'olya class $\sf LP$ is the only class of real entire functions closed under differentiation, 
but this result is still weaker than Wiman's conjecture W1. One down, six to go!

In the same work they also gave a detailed proof of P\'olya's conjecture P2 for entire functions 
of finite order~\cite[Theorem~2]{HW-1}. Their proof essentially followed {\AA}lander's very sketchy outline~\cite{A-C.R.}. 

\vspace{.1in}\noindent
In 1983, Hellerstein, Shen, and Williamson~\cite{HSW} proved P2 in full generality, describing entire functions $f$ with $f f' f''$ having real zeros only. 
In addition to two families of non-real entire functions with real zeros of all derivatives, they singled out two families
of non-real entire function of infinite order with $f f' f''$ having only real zeros:
\[
f(z) = C \exp\bigl[ e^{ i (az+b)} \bigr],  \quad  f(z) = C \exp\bigl[K( i (az+b) - e^{ i (az+b)}) \bigr],
\] 
where $C$ is a complex constant, $a$ and $b$ are real numbers, and $-\infty < K \le 1/4$.  For the functions of these 
two families, $f'''$ must have non-real zeros. This work significantly used the ideas of Levin and Ostrovskii.

\vspace{.1in}\noindent
In 1987, Craven, Czordas, and Smith~\cite{CCS} proved a result which is only
slightly weaker than the P\'olya--Wiman conjecture PW: if $f$ is of order less than $2$ and has finitely 
many non-real zeros, then all sufficiently high derivatives belong to $\sf LP$ (actually, in~\cite{P4}
P\'olya formulated Conjecture~PW in this weaker form). Then, in~\cite{CCS-1}, they extended this
result to entire functions of order $2$ and minimal type.

\vspace{.1in}\noindent
In 1990, Kim~\cite{K} completed their result and settled the P\'olya--Wiman conjecture PW
(and therefore, also the second Fourier--P\'olya conjecture FP2, which is equivalent to PW). 

\vspace{.1in}\noindent
In 2000, Ki and Kim~\cite{KK} settled the first Fourier--P\'olya conjecture FP1
and gave a simpler proof of the P\'olya--Wiman conjecture PW.

\vspace{.1in}\noindent
In 1989, Sheil-Small \cite{SS} proved the second Wiman's conjecture W2, and therefore
the first Wiman conjecture for entire functions of finite order. As a special case, this result resolves 
a long-standing puzzle: if $P$ is a real polynomial of degree $n >1$, then $P^2 + P'$ should have 
at least $p$ complex conjugate pairs of zeros, where $p$ is the least integer bigger than or equal to $n/2$.
For the first Wiman's conjecture W1, a gap between finite order and \eqref{1} remained.

\vspace{.1in}\noindent
In 2003, this gap was closed by Bergweiler, Eremenko, and Langley. They 
proved in~\cite{BEL} that, for functions of infinite order with all zeros real, $f''$ has infinitely many
non-real zeros, thus completing the proof of Wiman's conjecture W1. This work makes a substantial 
use of the ideas of Levin and Ostrovski~\cite{LO} and Sheil-Small \cite{SS}.

\vspace{.1in}
P\'olya's conjecture P3 turned out to be most difficult.
We mention two partial results. 

\vspace{.1in}\noindent
In 2002, Edwards and Hellerstein~\cite{EH} generalized
Sheil-Small's theorem. They introduced the class ${\sf W}_{2p}^*$ consisting of entire functions $f=pg$, were 
$g\in {\sf W}_{2p}$, and $p$ is a real polynomial, and showed that 
\vspace{.1in}

{\em If $f\in {\sf W}_{2p}^*$, then any derivative
$f^{(k)}$, $k\geq 2$, has at least $2p$ non-real zeros.}

\vspace{.1in}\noindent
In 1993,  Bergweiler and Fuchs \cite{BF} proved the following:
\vspace{.1in}

{\em If $f$ is a zero free real entire function
of infinite order, then $f''$ has infinitely many non-real zeros.}

\vspace{.1in}\noindent
In 2005, Langley \cite{L} complemented the result
of~\cite{BEL} and proved that 
\vspace{.1in}

{\em For any real entire function $f$ of infinite order with all zeros real, any derivative $f^{(k)}$, $k \ge 3$ has
infinitely many non-real zeros.}

\vspace{.1in}
This is probably the most difficult result in the area, building on the technique developed in~\cite{BEL}.

\vspace{.1in}\noindent
In 2006, the final step in the proof of P\'olya's conjecture P3 was made in \cite{BE}:

\vspace{.1in}
{\em For a real entire function of finite order which {\em does not} belong to ${\sf LP}^*$, 
the number of non-real zeros of $f^{(k)}$ tends to infinity as $k\to\infty$.}
\vspace{.1in}

By combining the aforementioned results,
one arrives to the following neat statement:
\vspace{.1in}

{\em For an arbitrary real entire function $f$ one of the following possibilities hold:
either zeros of $f^{(k)}$ are all real for sufficiently large $k$, 
or the number of non-real zeros of $f^{(k)}$ tends to infinity as $k\to\infty$.}

\vspace{.1in}
Thus the paper \cite{LO} had a lasting influence on the subsequent research
which eventually led to a complete proof of conjectures of P\'olya and Wiman.
Levin and Ostrovskii did not prove these conjectures themselves but they showed the correct
path towards the proofs, which was followed by many mathematicians who
eventually completed the project. Other results related to this range of questions can be found 
in~\cite[Section~VI.5]{GO}.

\vspace{.1in}

The machinery of Levin--Tsuji characteristics has found additional applications in the theory of meromorphic functions,
see, for example, Bergweiler--Eremenko~\cite{BE2}. Matsaev, Ostrovskii, and Sodin~\cite{MOS, MOS1} used 
Levin--Tsuji characteristics for estimates of logarithmic determinants and the Hilbert transform. 
Khabibullin~\cite{Khabib} applied them to find conditions for existence of certain subharmonic minorants.

\vspace{.1in}
We conclude this part with a question closely related to the P\'olya and Wiman problems, raised by Levin and Ostrovskii in~\cite{LO} and apparently still open.
Let $\sf P$ be the class of entire functions 
of the form $f(z)=e^{-\alpha z^2} g(z)$, where $\alpha \ge 0$, and $g$ is an entire function of
genus $1$, having no zeros in the upper half-plane $\mathbb H$. By Obreshkov's 
theorem~\cite[Theorem~VIII.4]{Levin}, $\sf P$ coincides with the locally uniform closure of $H$-polynomials,
that is, polynomials having no zeros in $\mathbb H$, so many properties of $H$-polynomials persist
for functions of class $\sf P$. In particular, it is also closed under differentiation, similarly to the Laguerre--P\'olya 
class $\sf LP$.  Levin and Ostrovskii stated the following conjecture:

\vspace{.1in}
\noindent LO. {\em Suppose $f$ and all its derivatives have no zeros in the upper half-plane. Then either $f\in {\sf P}$, 
or $f$ coincides with one of the exceptional non-real functions in P\'olya's conjecture} P2.

\vspace{.1in}\noindent
We should note that presumably due to an oversight, Levin and Ostrovskii did not mention these exceptional functions 
when asking this question.

\vspace{.2in}


\noindent
{\bf 3. Comb functions}
\vspace{.1in}

Consider the class $\frak F$ of real entire functions with real $\pm 1$-points. The functions of this class 
appear as Hill discriminants/Lyapunov functions of the second order periodic linear ODE.
In this capacity, they appear already in Lyapunov's work~\cite{Lyapunov1}.
Then this class was identified by Krein~\cite{Krein1953} in connection with the spectral theory of Krein's strings.
He also pointed out connections with other problems of analysis (functional Pell's equation, periodic continued 
fractions). This class also arises in Chebyshev-type extremal problems. 

Recall that, by Edrei's theorem, the functions 
of class $\frak F$ have finite exponential type. However, in applications, this non-trivial result is usually not used;
the needed growth estimates follow from the specifics of the question under consideration.

\paragraph*{Parametrization of class $\frak F$.}
In~\cite{MO}, Marchenko and Ostrovskii 
discovered a remarkable parametrization of class $\frak F$ by conformal mappings of the upper half-plane on a
comb domain, allowing them to find sets of free parameters that describe the spectrum of Hill operators,
and  periodic and anti-periodic inverse spectral problems for the Sturm-Liouville operators. Since then, 
the Marchenko--Ostrovskii parametrization has became indispensable in several areas of analysis.

\vspace{.1in}\noindent{\bf Theorem 1.}
{\em Class $\mathfrak F$ consists of all entire functions represented in the form
\begin{equation}
\label{comb}
u(z)=\cos\theta(z),
\end{equation}
where $\theta$ is a conformal mapping of the upper half-plane onto a comb domain of the form
\begin{equation}\label{comb2}
\mathfrak C = \{ z\colon \Ima z>0, \ p<{\rm Re\,  z <q}\,\}\backslash\bigcup_{p<k<q}[k\pi, k\pi+ ih_k],
\quad\mbox{with some}\  h_k\geq 0,
\end{equation}
where $-\infty\leq p < q\leq+\infty$, $\theta(\infty)=\infty$,
and $\theta$ is extended to the
lower half-plane by symmetry, $\theta(z) = \overline{\theta(\bar z)}$}.

\vspace{.1in}
In the case $-\infty<p<q<+\infty$, the set $\Lambda$ of $\pm 1$-points of the function $u$ is finite,
and $u$ is a polynomial of degree $q-p$. If only one of the numbers $p$ and $q$ is infinite, then
the set $\Lambda$ is bounded from one direction, either above or below. 
Let, for instance, $\min \Lambda>-\infty$, and set $\mu = \min\Lambda$. 
Then the function $u((z-\mu)^2)$ has real $\pm 1$-points, and therefore, has a finite
exponential type. So, in this case, $u$ has order $1/2$ and finite type.

The idea of the proof of Theorem 1 is elegant and not very difficult. Consider the case when the set
$\Lambda$ is unbounded both from below and from above. By Edrei's
theorem, the functions $u\pm 1$ belong to the Laguerre-P\'olya class.
Hence, by the Laguerre theorem, $u'$ has only real zeros, which interlace both with $+1$- and with $-1$-points of $u$.
So we can enumerate solutions to the equations $u^2=1$ and $u'=0$ so that
\[
\ldots < a_k \le  c_k \le b_k < a_{k+1} \le c_{k+1} \le b_{k+1} < \ldots\,,
\]
where $u(a_k)=u(b_k)=(-1)^k$, and $u'(c_k)=0$. Consider the meromorphic function
\[
[(\arccos u)']^2 = \frac{(u')^2}{1-u^2}\,.
\]
Its zeros and poles coincide with the ones of the convergent infinite product
\[
\prod_k \frac{(1-z/c_k)^2}{(1-z/a_k)(1-z/b_k)}\,.
\]
Applying growth considerations, it is not difficult to see that
\[
\frac{(u'(z))^2}{1-u(z)^2} = C\prod_k \frac{(1-z/c_k)^2}{(1-z/a_k)(1-z/b_k)}\,,
\]
with some $C>0$, whence
\[
(\arccos u(z))' 
= D \prod_k \frac{1-z/c_k}{\sqrt{(1-z/a_k)(1-z/b_k)}}\,,
\]
with $D>0$, and an appropriate choice of the branches of the square roots. By the Schwarz--Christoffel theorem,
the functions
\[
\theta_n (z) = \int_0^z \prod_{|k|\le n} \frac{(1-\zeta/c_k)}{\sqrt{(1-\zeta/a_k)(1-\zeta/b_k)}}\,  d\zeta
\]
are conformal maps of the upper half-plane onto comb domains, and their
$n\to\infty$ limit $\theta$ is also univalent in the upper-half plane and maps it onto a comb domain.

\paragraph*{Spectrum of Hill's operator.}
Consider Hill's equation
\begin{equation}\label{eq:Hill}
Ly = -y'' + q(x)y = \lambda y 
\end{equation}
with a real periodic potential $q\in L^2(0, \pi)$,  $q(x)=q(x+\pi)$. 
The main problem is to determine {\em stability intervals} that consist of values of $\lambda$ for which
all solutions to~\eqref{eq:Hill} are uniformly bounded on the real axis. It is known that the closure of the union of stability intervals coincides with the spectrum of Hill's operator $L$ acting in $L^2(\R)$ (see, for instance, Glazman~\cite[Section~56]{Glazman} or Luki\'c~\cite[Section~11.16]{Luk}), so, in the spectral theory language,
the problem is to describe the spectra of Hill's operators.

Let us start with a brief reminder of a fragment of Floquet's classical theory.
Let $c(x,\lambda)$ and $s(x,\lambda)$ be the fundamental solutions to~\eqref{eq:Hill} with initial conditions
$$s(0,\lambda)=c'(0,\lambda)=0,\quad c(0,\lambda)=s'(0,\lambda)=1.$$
Then $c$ and $s$ are entire functions of $\lambda$. For large complex $\lambda$, they are close
to $\cos (x\sqrt{\lambda})$ and $\sin (x\sqrt{\lambda})/\sqrt{\lambda}$,
 respectively.
In particular, both of them have order $1/2$, mean  type, see, for instance, Levitan--Sargsjan~\cite[Chapter~1, \S2]{LS}.
For an arbitrary solution $y$ to~\eqref{eq:Hill}, we have
\[
y(x) = y(0)c(x, \lambda) + y'(0)s(x, \lambda),
\]
whence
\begin{equation}\label{eq:monodromy}
\left(
\begin{matrix}
y(\pi) \\
y'(\pi)
\end{matrix}
\right)
= T(\lambda)
\left(
\begin{matrix}
y(0) \\
y'(0)
\end{matrix}
\right),
\end{equation}
where
$$T(\lambda)=\left(\begin{array}{cc}c(\pi,\lambda)& s(\pi,\lambda)\\
c'(\pi,\lambda)& s'(\pi,\lambda)\end{array}\right)$$
is the monodromy matrix of Hill's equation. Thus, the problem of determining stability intervals boils down to 
finding the eigenvalues of the matrix $T(\lambda)$, which we denote by $\rho_1$, $\rho_2$.

The determinant of $T(\lambda)$ identically equals $1$ (since it is
the Wronskian 
of the functions $c$ and $s$, and equation (\ref{eq:Hill})
contains no term with $y'$). So $\rho_1 \cdot \rho_2 =1$, and
\[
\rho_1 + \rho_2 = \operatorname{trace} T(\lambda) = c(\pi, \lambda) + s'(\pi, \lambda).
\]
We consider $1/2$
of this trace, which is called the {\em Lyapunov function} (or the {\em Hill discriminant}),
$$u(\lambda)=(c(\pi,\lambda)+s'(\pi,\lambda))/2.$$
This is a real entire function of order $1/2$, mean type. The key observation, known already to
Lyapunov~\cite{Lyapunov1}, is that all $\pm1$-points of
$u$ are real.
Indeed, if $u(\lambda)=1$, then,
$\rho_1=\rho_2 = 1$, and, by~\eqref{eq:monodromy}, Hill's equation has a periodic solution,
that is, $\lambda$ belongs to the spectrum of the periodic boundary values problem for Hill's operator.
Similarly,  if  $u(\lambda)=-1$, then
$\rho_1=\rho_2 = -1$, and $\lambda$ belongs to the spectrum of the anti-periodic boundary 
values problem. Both boundary values problems are self-adjoint, so, in both cases, $\lambda$ 
must be real.  We also note that, $u(\lambda)\to+\infty$ as $\lambda\to-\infty$;
this follows from an asymptotic analysis of solutions $c$ and $s$, see, for instance,
\cite[Ch. I, \S2]{LS}. So the set of $\pm1$-points of $u$ is bounded from below. 

Then, by a straightforward inspection, the stability set contains  intervals of the real axis, where $|u|<1$, 
in which case, $\rho_1\ne \rho_2$, and both are unimodular, 
while the instability set contains intervals where $|u|>1$, in which case $\rho_1\ne \rho_2$, and
both are real. At the points $\lambda$, where $u(\lambda)=\pm 1$, we have $\rho_1 = \rho_2$, 
and the answer depends on whether the derivative $du/d\lambda$ vanishes at this point or not.  
If $ du/d\lambda \ne 0$, then $\lambda$ belongs to the instability set. 
If $ du/ d\lambda=0$ (i.e., two neighbouring
stability intervals stick together), then $\lambda$ belongs to the stability set. All this was known by the end of the
19th century, in particular, to Lyapunov. The absent piece was a
parametric description of Lyapunov functions corresponding to
Hill equations. 

\vspace{.1in}
Representation of the Lyapunov function in the form \eqref{comb} combined with accurate estimates of
the solutions $c$ and $s$ allowed Marchenko and Ostrovskii to use the well-developed theory of distortion under 
conformal maps and to obtain complete characterization of Lyapunov's functions, and hence a characterization 
of the spectra, of Hill's operators with real potentials $q$ in the Sobolev space
$\widetilde{W}_2^n[0, \pi]$ of $\pi$-periodic functions, which consists of
functions $v$ on $[0, \pi]$ satisfying
$$\sum_{k=0}^n\int_0^\pi |v^{(k)}(x)|^2\, dx<\infty$$
and continued $\pi$-periodically on $\R$.
\vspace{.1in}

\noindent
{\bf Theorem 2.} {\em If a real potential $q$ belongs
to $\widetilde{W}_2^n[0, \pi]$ with $n \ge 0$, 
then the Lyapunov function $u$ has the form \eqref{comb} with $p>-\infty$, 
$\theta(x(1+i))\sim\pi\sqrt{x}$, $x\to+\infty$, and
\begin{equation}\label{eq:slits}
\sum_{k=p+1}^\infty(k^{n+1}h_k)^2<\infty.
\end{equation}
Conversely, for every conformal map $\theta$ of the upper half-plane onto
a comb domain that satisfies these conditions, the function $u$, as defined in
\eqref{comb}, is the Lyapunov function of Hill's operator with real potential
in $\widetilde{W}_2^n[0, \pi]$.}
\vspace{.1in}

Theorem 2 has several remarkable consequences in spectral theory, see~\cite[\S4, Chapter~3]{M}. 
Probably, the most outstanding one is the identification of the moduli space of Hill's operators
with potentials in the Sobolev space 
(that is, independent spectral data, which allow one to uniquely recover the potential)
with the following set of data: 

\smallskip\noindent (i)
a comb domain $\mathfrak C$ with slits satisfying~\eqref{eq:slits};

\smallskip\noindent (ii)
a marked point on each slit (if a slit degenerates, i.e., $h_k=0$, then the point $\pi k$ is marked); 

\smallskip\noindent (iii)
a sign attached to each marked point, not lying on the base of the slit.

\vspace{.1in}
Let us sketch the main idea behind this identification, see~\cite{MO} and \cite[Chapter 3 \S4]{M} for the details
(cf. Stankevich~\cite{Stankevich}). 
Suppose that we know the comb-domain $\mathfrak C$, that is, the Lyapunov function $u$ of Hill's equation~\eqref{eq:Hill}. 
Zero sets of the functions $s$ and $s'$ are eigenvalues of Hill's operator with boundary conditions
$y(0)=y(\pi)=0$, and $y'(0)=1$, $y(\pi)=0$. Classical methods of the inverse spectral theory recover the 
potential $q$ by these two spectra. Suppose that we know zeros of $s$. At these points,
by unimodularity of the matrix $T$, $c(\lambda)s'(\lambda)=1$, and therefore,
$|u(\lambda)| = |c(\lambda) + s'(\lambda)|/2 \ge 1$. Thus, zeros of $s$ lie on the slits of $\mathfrak C$, 
and each slit contains exactly one zero of $s$. The $\theta$-images of zeros of $s$ are the marked points from item (ii).

Since $s$ is an entire
function of order $1/2$, we can recover $s$ up to a multiplicative constant as an infinite product of genus zero.
The constant can be found from the asymptotics $s(\lambda)\sim \sin(\pi\sqrt{\lambda})/\sqrt{\lambda}$, 
$\lambda\to - \infty$. Consider the function $v=(c-s')/2$. Having $v$, we recover $s'=(u+v)/2$, and hence, zeros
of $s'$. Thus, the problem boils down to recovery of $v$. 

The next step is to expand the meromorphic function $v/s$ into a series of simple fractions.
Such an expansion yields that 
$v$ can be recovered by knowledge of it at zeros of $s$. Using again that $cs'=1$ at zeros of $s$,
we see that $u^2-v^2=1$ therein, and therefore,
$ v(\lambda) = [{\rm sgn\,}v(\lambda)]\, \sqrt{u^2(\lambda)-1} $,
provided that $s(\lambda)=0$. Thus, knowing signs of $v(\lambda)$ at zeros of $s$, we can recover 
the values of $v$ at these points, and then $v$ itself.

 \vspace{.1in}
The Marchenko--Ostrovskii theory generated a very large body of literature which is impossible
to survey here. So we mention only a few developments stemming from~\cite{MO}.

\paragraph*{Related spectral problems.}
Korotyaev~ \cite{Kor} extended Theorem 2 to a wider class of potentials, 
real distributions of the form $\{ q=v'\colon v\in L^2(0,\pi)\}$ continued with period $\pi$ on $\R$.

Tkachenko \cite{T}
obtained  a parametrization of spectra of Hill's operators whose 
potential $q$ is not necessarily real, so these operators are not
self-adjoint. In this case comb representation is not available, and
one has to work directly with the Riemann surface of $u^{-1}$ spread over
the plane.

Already in 1953 Krein~\cite{Krein1953} studied spectral properties of solutions to the integral equation
\[
y'(x) + \lambda \int_0^x y(t)\,  d M(t) = \operatorname{const},
\]
with a non-decreasing function $M$, $M(x+1)-M(x)=\operatorname{const}$. This equation describes the amplitude
$y(x)$ of small vibrations of a string with periodic distribution of mass. 
Krein stated that an entire function $u$ normalized by $u(0)=1$
is Lyapunov's function for a periodic string equation if and
only it has only non-negative $\pm 1$-points.
Actually, Krein imposed an extra condition that $u$ is represented by convergent Hadamard product of
genus zero with positive zeros, but in view of Edrei's theorem, this condition is redundant. Combining Krein's theorem
with Theorem 1, one gets a description of spectra of periodic Krein's string.

Mikhailova~\cite{Mikhailova} used the Marchenko--Ostrovskii map to parameterise 
monodromy matrices of two-dimensional canonical systems with a periodic Hamiltonian, and gave a constructive 
procedure of recovery of the monodromy matrix.

\paragraph*{Finite-band potentials and their closure.}
It is a natural hope to extend the Marchenko--Ostrovskii theory from periodic potentials to
more general classes of functions
invariant with respect to action of the real axis by translations.
Usual suspects
are almost-periodic and, more generally, random ergodic potentials.
There are only a few classes of such potentials for which 
this dream has been fully realized.

Given a real potential $q$, we denote by $L$ the Schr\"odinger self-adjoint operator
\[
L = - \frac{ d^2}{ dx^2} + q
\]
acting in $L^2(\R)$. By $\sigma (L)$ we denote its spectrum.

First, we consider the class of {\em finite-band} potentials which we define shortly. 
A concise introduction to their theory can be found in Moser lectures~\cite{Moser1}. 
See also Akhiezer~\cite{Akhiezer-Inverse} and Dubrovin, Matveev, and Novikov~\cite[Chapter~2]{DMN}.
An interesting discussion of the history is in the article by
Matveev~\cite{Matveev}.
These potentials have a finite-band spectrum 
\begin{equation}\label{eq:sigma}
\sigma (L) = [\lambda_0, \lambda_1] \cup [\lambda_2, \lambda_3] \cup \ldots [\lambda_{2N}, \infty), \quad \lambda_0<\lambda_1<\lambda_2 < \ldots < \lambda_{2N},
\end{equation}
and are singled out by a resolvent condition. The resolvent $R(\lambda) = (L-\lambda)^{-1}$.
is an integral operator with the kernel $G(x, y, \lambda)$ called Green's function of $L$. 
Its diagonal $G(x, x, \lambda)$ belongs to the Nevanlinna class~$\mathsf N$ of functions $f$, analytic in 
$\{ \Ima\lambda\ne 0\}$, and satisfying $f(\bar\lambda) = \overline{f(\lambda)}$, and 
$\Ima f(\lambda)/\Ima\lambda>0$. 

The potential with spectrum~\eqref{eq:sigma} is called finite-band if its Green function $G(x, x, \lambda)$
takes purely imaginary limiting values on the interiors of the bands. Then it is not difficult to show that
$G(x, x, \lambda)$  can be continued to a meromorphic function of $\lambda$ on the 
hyperelliptic Riemann surface $X$ of the function
\[
\Bigl[\, \prod_{0\le j \le 2N} (\lambda-\lambda_j) \,\Bigr]^{1/2}.
\]
A remarkable consequence is that, for finite-band potentials, the direct and inverse problems 
of spectral theory can be explicitly solved by classical methods of the theory of Riemann surfaces.
The full set of spectral data 
of finite-band  potentials is realized as divisor on $X$,
which consists of $N$ points chosen so that each closed gap $[\lambda_{2j-1}, \lambda_{2j}]$
has a point lying over it. The potential $q$ is recovered using theta-functions (or hyperelliptic integrals), 
and is quasi-periodic (recall that a function $f$ is called {\em quasi-periodic} if 
$f(x) = F(x, \ldots , x)$, with a continuous function $F(x_1, \ldots , x_n)$ that is periodic
in $x_1, \ldots , x_n$. Quasi-periodic functions are
uniformly almost-periodic).

Alternatively, one can define the spectral data using comb-domains.
In this case, a function $\theta$ maps the upper half-plane onto a quater-plane with $N$ slits having
{\em arbitrary} bases.  
Then, as in the Marchenko--Ostrovskii theory, the complete set of spectral invariants of the operator $L$ is given by 
location and lengths of slits,  collection of $N$ points on these slits, and collection of signs attached to each point 
unless the point lies on the base of a slit.  

If the bases of the slits lie in a subset of an arithmetic progression, then the finite-band potential is periodic and
we arrive at a special case of the Marchenko--Ostrovskii theory with Lyapunov function $u$ having  finitely many simple 
$\pm 1$-points, the rest of them having multiplicity two.
Using Theorem~2, Marchenko and Ostrovskii showed that such 
finite-band periodic potentials are dense in the $\widetilde{W}_2^n[0, \pi]$-metric in the space of all real periodic potentials~\cite{MO2}, \cite[Theorem~3.4.3]{M}.

Comb domains also occur for some classes of almost-periodic potentials with Cantor-type spectrum
when  the bases of the slits are everywhere dense, see Pastur and Tkachenko~\cite{PT}, and 
Sodin and Yuditskii~\cite{SY1} (the latter work contained a gap, which was fixed in
Gesztesy and Yuditskii~\cite{GY}). In these works, the authors obtained a complete description of
independent spectral data needed for the recovery of the potential.

\paragraph*{Almost-periodic and random ergodic potentials.}
For general almost periodic and random ergodic potentials the classical Floquet theory ceases to work, and the Lyapunov
function $u$ does not exist. Johnson and Moser~\cite{JM} found an extension of the Floquet theory for almost-periodic 
potentials, and Kotani~\cite{Kotani1, Kotani2, Kotani3} extended it to random ergodic potentials defined as follows.
Let $\mathbb R$ (regarded as an additive group) act on a probability space $(\Omega, \mathcal F, \mathbb P)$ 
by measure preserving transformations $T_x$, $x\in\mathbb R$, let the action be ergodic, and  
let $q\colon \Omega\to\mathbb R$ be a random variable with a finite second moment. 
Then one considers Hill's operators, whose potentials are the random functions $x\mapsto q(T_x\, \omega)$.

We will briefly explain the approach developed by Johnson--Moser and Kotani. 
If the potential $q$ is either uniformly almost periodic, or, more generally, random ergodic potential $q$,
then the Green's function $x\mapsto G(x, x; \lambda)$, as well as its inverse $x\mapsto 1/G(x, x; \lambda)$, 
are also almost-periodic/random ergodic. Set
\begin{equation}\label{eq:w}
w(\lambda) \stackrel{\rm def}= - \frac12\, \lim_{b-a\to\infty} \frac1{b-a}\, \int_a^b \frac{ dx}{G(x, x; \lambda)}.  
\end{equation}
In the random case the limit exists almost surely, and, by ergodicity, coincides with
$ - \frac12\, \mathbb E \bigl[ G^{-1}(x, x; \lambda) \bigr]$. The function $w$ is holomorphic in 
$\{\Ima\lambda \ne 0\}$ and satisfies
\[
w, \, \frac{ dw}{ d\lambda}, \, - iw \in \mathsf N,
\]
where ${\mathsf N}$ is the Nevanlinna class introduced above.
This yields that $w$ is univalent in 
the upper half-plane, and the image $w(\mathbb C_+)$  lies in the second quadrant. Generally speaking, 
$w(\mathbb C_+)$ is not a comb domain, but still, if $\lambda\in w(\mathbb C_+)$, 
then, for any $\tau>0$, $\lambda-\tau\in \omega(\mathbb C_+)$,  as well.

To describe the boundary of the image $w(\mathbb C_+)$, one needs to know the boundary values of
$w$ on the real axis. They are equal  to $-\gamma +  i\pi k$, where
$\gamma$ is {\em the Lyapunov exponent}, and $k$ is
{\em the integrated density of states}. 
The Lyapunov exponent $\gamma(\lambda)$ measures the exponential
growth of solutions
\begin{equation}\label{eq:Lyapunov}
\gamma (\lambda) = \lim_{x\to\pm\infty} \frac1{|x|}\,
\log\| T(x, \lambda) \|, 
\end{equation}
where $T$ is the fundamental matrix,
\[
T(x, \lambda) = 
\left(
\begin{matrix}
c(x, \lambda) & s(x, \lambda) \\
s'(x, \lambda) & c'(x, \lambda)
\end{matrix}
\right).
\]
The integrated density of states is defined as
\[
k(\lambda) = 
\lim_{b-a\to \infty}\, \frac{\nu (a, b; \lambda)}{b-a}\,,
\]
where $\nu (a, b; \lambda)$ is the number of the eigenvalues of the Dirichlet eigenvalue problem 
on the interval $[a, b]$, which are less or equal to $\lambda$. If the potential $q$ is periodic with period $\pi$, then 
$w=  i\theta$ where $\theta$ is the Marchenko--Ostrovskii conformal map.
In the case of random Jacobi matrices, the conformal map $e^w$ was thoroughly studied by Hur and Remling~\cite{HR}.

The Lyapunov exponent and integrated density of states contain important spectral information 
about Schr\"odinger operators with almost-periodic and random ergodic potentials, see the Pastur--Figotin 
treatise~\cite[Chapter~V]{FP}. For instance,
the spectrum of $L$ coincides with the closed support of the measure $ dk$, while the
support of the absolutely continuous spectrum of $L$ coincides with the essential closure of the set 
$\{\gamma=0\}$ with respect to the Lebesgue measure (that is, with the set of points $x\in\R$ such that 
any neighbourhood of $x$ has an intersection of positive Lebesgue measure with $\{\gamma=0\}$). It is 
worth mentioning that $\gamma$ is a non-negative subharmonic function in $\mathbb C$ with the Riesz measure 
$ dk$, and that $ dw/ d\lambda$ is the Hilbert transform of $ dk$.

\paragraph*{Chebyhsev-type extremal problems.}
Comb functions of the form~\eqref{comb}, both polynomial and transcendental, naturally appear as extremal functions 
in various Chebyshev-type problems. They were introduced by Akhizer and Levin in~\cite{AL}. Using comb functions, 
Akhiezer and Levin obtained far reaching generalizations of Bernstein's inequality for entire functions of exponential type. 
Note that Akhiezer and Levin considered a more general class of comb domains, for which bases of the slit 
do not have to be a subset of an arithmetic progression. In that case, the functions of the form~\eqref{comb} are analytic 
in the plane with cuts along $\theta^{-1}\R$. These functions also play an important role in Levin's theory of subharmonic 
majorants~\cite{Levin-M}.

Employing comb functions, Eremenko \cite{Eremenko-DAN} found the upper envelope $M(x)$, $x\in\R$, of absolute values 
of entire functions $f$ of a given exponential type, such that 
\[
|f(x)| \le
\begin{cases}
A, & x<0; \\
B, & x\ge 0.
\end{cases}
\] 
Note that, for complex values $z$, the exact value of the majorant $M(z)$ is not known (Hayman's question).

A rather general result connecting Chebyshev-type extremal problems for entire functions with comb domains 
can be found in the Sodin--Yuditskii survey paper~\cite[Theorem~7.5]{SYu}.

\paragraph*{Pell's equation.}
Chebyshev and his pupils reduced many polynomial extremal problems to polynomial analogues of Pell's equation,
\begin{equation}\label{eq:Pell}
X^2 + TY^2 = 1,
\end{equation}
see~\cite{SYu}.  Krein~\cite{Krein1953} found a relation between entire functions, which appear in the spectral
theory of Krein's string, and Pell's equation for entire functions. Akhiezer~\cite{Akhiezer-Pell} demonstrated an
approach to construction of finite-band potentials via functional Pell's equations.  
So it is not too surprising that comb polynomials and  entire functions are related to Pell's equation in entire functions.
This relation was clarified by Marchenko and Ostrovskii in~\cite[Theorem~6.1]{MO}.

Denote by $\mathsf C$ the class of real entire functions of Cartwright class (defined in the beginning of Section~1)
with real zeros. 
\vspace{.1in}

\noindent
{\bf Theorem 3.}
{\em Let $T\in \mathsf C$, and $T(0)>0$. Then Pell's equation~\eqref{eq:Pell} has a solution $X, Y\in\mathsf C$ iff
and only if
there exists an entire function $\Phi\in\mathsf C$ and a constant $b$ such that the function
\[
\theta (z) = \int_0^z \frac{\Phi (\zeta)}{\sqrt {T(\zeta)}}\,  d\zeta + b
\]
maps conformally the upper half plane onto a comb domain $\mathfrak C$.

Furthermore,  a general form of solutions $X, Y\in\mathsf C$ is
\[
X(z) = \cos\theta (z), \quad Y(z) = \frac{\sin\theta (z)}{\sqrt{T(z)}}\,,
\]
where $\theta$ maps the upper half-plane onto a comb domain $\mathfrak C$.
}

\vspace{.1in}
In the polynomial case, the solvability of Pell's equation is equivalent to existence of a special
factorization $1=(X+ i\sqrt{T}Y)(X- i\sqrt{T}Y)$ of the function  identically equal to $1$
on the hyperelliptic Riemann surface of $\sqrt{T}$. In the transcendental case, we arrive at
the same factorization on a hyperelliptic Riemann surface of infinite genus.  

\paragraph*{Entire functions of Krein's class.}
In~\cite{O}, Ostrovskii found a parametric description of Krein's class $\sf K$ (we defined it 
in the very beginning of Section~1). Note that, since functions from Krein's class belong to the Cartwright class
and are real, they have representation
\[
f(z) = Cz^\ell \lim_{R\to\infty} \prod_{\Lambda\cap [-R, R]} \Bigl( 1 - \frac{z}\lambda \Bigr),
\]
where $\ell=0$ or $1$ and $\Lambda$ is the zero set of $f$. Hence, functions from Krein's class 
are defined by their zeros up to a constant factor and, for description of the class $\sf K$, it suffices to describe 
the sets $E\subset\mathbb R$ that are zero sets of functions from $\sf K$.

Let us say that a sequence $E$ of real numbers $\lambda_k$, $-\infty\leq p<k<q\leq\infty$ satisfies an $R$-condition 
if there is a function $\theta$ as in Theorem 1, for which $\lambda_k\in[a_k,b_k]$ where $[a_k,b_k]$ is
the $\theta$-preimage of the $k$-th ``tooth'' $[\pi k, \pi k+ih_k]$.

\vspace{.1in}
\noindent
{\bf Theorem 4.} {\em A set $E$ coincides with the zero set of some function
of class $\sf K$ if and only if it satisfies the $R$-condition.}

\vspace{.1in}
The importance of this theorem lies in the fact that, as shown by Krein in~\cite{Krein2},
functions from class $\sf K$ describe the entries 
of entire Nevanlinna matrices. A unimodular entire matrix-function
\[
\left(
\begin{matrix}
A(z) & B(z) \\
C(z) & D(z)
\end{matrix}
\right)
\]
is called a Nevanlinna matrix if its entries are real entire functions, and, for all real values $t$, the function
\[
z \mapsto \frac{A(z)t+B(z)}{C(z)t+D(z)}
\]
has positive imaginary part in the upper half-plane.
These matrices play a fundamental role in the description of solutions of many 
problems  in analysis, which include the Hamburger moment problem, the spectral theory of Schr\"odinger-type linear differential equations, Krein's problem on continuation of positive definite functions from an interval, and de Branges theory. The monodromy matrices $T$, which we dealt with above, are also Nevanlinna matrices.

In the same paper~\cite{O}, Ostrovskii  
showed that $\mathfrak F \subsetneq \sf K$ (recall that the class $\mathfrak F$ consists of real entire functions with real 
$\pm 1$-points), and that every function of the Krein's class can be represented as
a sum of two comb functions from $\mathfrak F$. The latter result was improved by Katznelson~\cite{Katz}.
He showed that {\em every real  entire function $f$ of Cartwright class can be represented as a sum of
two comb-entire functions: $f=u_1+u_2$, $u_1,  u_2\in\mathfrak F$}. Moreover, the exponential types of the
functions $u_1$, $u_2$ coincide with that of $f$.

\paragraph*{Vinberg combs. P\'olya functions.}
Another useful comb representation of entire functions was introduced by Vinberg~\cite{V}.
Given $-\infty \le p < q \le +\infty$, let $(c_k)_{p<k<q}$ be a sequence of real numbers.
Consider the region
\begin{equation}\label{vcomb}
\mathfrak V=\{ z\colon \{\pi  p < \Ima z < \pi q\}\,  \big\backslash
\bigcup_{p < k < q} \bigl\{  i\pi k+t\colon -\infty<t<\log|c_k| \bigr\}.
\end{equation}
Such regions are called Vinberg's combs. Let $\theta$ be a conformal
map of the upper half-plane onto a Vinberg's comb normalized by 
$\displaystyle \lim_{y\to +\infty}\theta( iy)=+\infty$. Then by the Symmetry Principle the function
\begin{equation}\label{vcomb2}
f=\exp\theta
\end{equation}
extends to an entire function. This entire function is real on the real line,
has critical values $\pm c_k$, and possibly one asymptotic value $0$.
\vspace{.1in}

\noindent
{\bf Theorem 5.} {\em The class of functions (\ref{vcomb2}) coincides
with the Laguarre--P\'olya class.}
\vspace{.1in}

As a corollary Vinberg obtains a parametrization of
the Laguerre-P\'olya class by
sequences of critical values. Comb functions of Marchenko and Ostrovskii
correspond to the subset
with $|c_k|\geq 1$, and for them we have both
representations  (\ref{comb}) and (\ref{vcomb2}). Using his comb representation, Vinberg obtained 
a simple purely topological proof of the results of
MacLane~\cite{Mc}
on real entire functions with real critical points.
See the survey~\cite{EY} for details and
examples.

\vspace{.1in}
In~\cite[Chapter~I]{deBr}, de Branges noticed that the theory
of entire functions
of the Laguerre-P\'olya class can be built starting with the following
observation:
{\em if $p$ is a polynomial without zeros in the upper half-plane, then, for any $x\in\R$, 
the function $y\mapsto |p(x+ iy)|$ increases with $y$}, which, in turn, yields that $p'/p$ 
has negative imaginary part in the upper half-plane.

We take this property as a starting point and say that {\em a function $f$ analytic in the upper 
half-plane $\mathbb C_+$ belongs to the  P\'olya class $\mathfrak P$ if the function $y\mapsto |f(x+ iy)|$
does not decrease, as $y$ increases}. Clearly, $\sf{LP}\subsetneq \sf{P} \subsetneq \mathfrak{P}$
(the class $\sf P$ was introduced in the very end of Section~2).
In a sense, these classes are not too far from each other. If  $f\in\mathfrak P$, 
then the ``mirror continuation'' 
of $\log |f|$ defined as
\[
v(z) = 
\begin{cases}
\log |f(z)|, & \Ima z>0, \\
\log |f(x+ i0)|,  & x\in\R, \\
\log |f(\bar z)|, & \Ima z<0,
\end{cases}
\] 
is subharmonic in $\C$. Informally speaking, $v$ belongs to a subharmonic counterpart of the 
Laguerre--P\'olya class {\sf LP}.

The functions of P\'olya class were independently introduced and
studied by Kargaev and Korotyaev~\cite{KarKor} and Weitz~\cite{Weitz}.
One of their results was a parametrization of class $\mathfrak P$ by 
conformal maps onto domains $\Omega$ such that {\em if $\zeta\in\Omega$, then, for any 
$t>0$, $\zeta+t\in \Omega$}.

\vspace{.1in}
\noindent
{\bf Theorem 6.} {\em 
$f\in\mathfrak P$ if and only if $\log f$ is a conformal map from the upper half plane onto a domain 
$\Omega$ such that $\lim_{y\to +\infty}{\rm Re} \log f( iy)=+\infty$.
}

\vspace{.1in}
The key step in proving Theorem 6 is the following observation:
{\em $f\in \mathfrak P$ if and only if either $\Ima\, f'/f <0$ in the upper half-plane, or
$f(z)=e^{az+b}$ with $a\in\R$}.

Theorem 6 contains as a special case Vinberg's Theorem~5. It also contains the
results
of Levin~\cite{Levin-M}. In connection with his theory of subharmonic majorants, Levin
considered the case when $\mathfrak P$ lies in the right half-plane. We also note that if $w$ 
is a function introduced by Johnson--Moser and Kotani and defined in~\eqref{eq:w}, 
then $w=-\log f$, with $f\in \mathfrak P$.

\paragraph*{Acknolwedgments.}
We are grateful to Jim Langley, Sasha Sodin, and Peter Yuditskii,
who read a draft of this paper and provided many comments that we have taken into account.

\vspace{.2in}

\noindent A.E.: Department of Mathematics, Purdue University, West Lafayette, IN 47907, USA \newline
{\tt eremenko@purdue.edu}

\vspace{.1in}
\noindent M.S.: School of Mathematics, Tel Aviv University, Tel Aviv, Israel \newline
{\tt sodin@tauex.tau.ac.il}

\end{document}